
\documentstyle{amsppt}

\def\a{\kern+.6ex\lower.42ex\hbox{$\scriptstyle \iota$}\kern-1.20ex a}
\def\e{\kern+.5ex\lower.42ex\hbox{$\scriptstyle \iota$}\kern-1.10ex e}
\def\z{{\accent95 z}}
\define\Jac{\operatorname{Jac}}
\define\Aut{\operatorname{Aut}k^2}

\topmatter
\title Injectivity on one line \endtitle
\author by Janusz Gwo\'zdziewicz \endauthor
\address
Department of Mathematics,\newline
Technical University,\newline
Al.~1000\,LPP\,7, 25--314~Kielce, \newline
Poland
\endaddress
\abstract
Let $k$ be an algebraically closed field of characteristic zero.
Let $H:k^2\to k^2$ be a polynomial mapping such that the
Jacobian $\Jac H$ is a non-zero constant.  In this note we
prove, that if there is a line $l \subset k^2$ such that
$H|_l:l\to k^2$ is an injection, then $H$ is a polynomial
automorphism.
\endabstract
\endtopmatter

\document
\head 		1. Main result
\endhead

Let $k$ be an algebraically closed field of characteristic zero.
Put $k^{*}=k\setminus \{0\}$.
By $\Aut$ we denote the group of polynomial automorphisms i.e. all
mappings $H=(f,g):k^2 \to k^2$, $f$, $g\in k[x,y]$ for which
there exists an inverse polynomial mapping.
Remind that for $H\in \Aut$ the Jacobian $\Jac H$
is a non-zero constant.
The famous Keller conjecture states that any polynomial mapping
which has the non-zero constant Jacobian is a polynomial automorphism.

\proclaim{Theorem~1.1}
Let $H:k^2\to k^2$ be a polynomial mapping such that $\Jac H\in k^{*}$.
If there exists a line $l\subset k^2$ such that $H|_l:l\to k^2$ is
injective then $H$ is a polynomial automorphism.
\endproclaim

The proof of the above theorem is given in Section~2 of this
note. Let us note here that a weaker result (injectivity on
three lines) has been proved recently in \cite 4.

\head 		2. Proof of the theorem
\endhead

The proof of our result is based on the Abhyankar--Moh theorem
and on some properties of Newton's polygon which we quote below.

If $f$ is a polynomial in $k[x,y]$ then $S_f$ denotes the
support of $f$, that is, $S_f$ is the set of integer points $(i,j)$ such
that the monomial $x^iy^j$ appears in $f$ with a non-zero
coefficient. We denote by $N_f$ the convex hull (in the real
space $\Cal R^2$) of $S_f \cup \{(0,0)\}$. The set $N_f$ is
called (see \cite 2) {\it Newton's polygon \/} of $f$.

\proclaim{Theorem 2.1}(\cite 2, \cite{3, theorem~3.4.})
Let $f$, $g$ be polynomials in $k[x,y]$. Assume that $\Jac (f,g)$ is
a non-zero constant and $\deg f>1$, $\deg g>1$.
Then the polygons $N_f$ and $N_g$ are similar, that is
$N_g = \frac{\deg g}{\deg f}N_f$.
\endproclaim

The following lemma is easy to check, so we omit the proof.

\proclaim{Lemma 2.2}
Let $f$, $g$ be polynomials in $k[x,y]$ with $\Jac (f,g)\in k^{*}$.
If $\deg f\le 1$ or $\deg g\le 1$ then the mapping $(f,g)$ is
a polynomial automorphism.
\endproclaim

\proclaim{Lemma 2.3}
Let $H=(f,g)$ be a polynomial mapping such that $\Jac H$ is a
non-zero constant.  If $H(x,0)=(x,0)$  then $H\in \Aut$.
\endproclaim

\demo{Proof}
First suppose that $\deg f>1$, $\deg g>1$. We have
$$ \align
   f(x,0)=x  \tag 1 \\
   g(x,0)=0. \tag 2
\endalign $$
{}From (1) the point $(1,0)$ belongs to the polygon $N_f$, so by
theorem~2.1 $(\frac{\deg g}{\deg f},0)\in N_g$.
This means that the polynomial $g$ contains some monomials of
the form $x^i$, $i>0$ with non-zero coefficients but this is a
cotradiction with (2).

We have $\deg f\le 1$ or $\deg g\le 1$ and by lemma~2.2 $H$ is a
polynomial automorphism.
\enddemo

\demo{Proof of theorem~1.1}
Without loss of generality we may assume that the line $l$ has
an equation $y=0$. Otherwiese we replace $H$ by $H\circ L$ where
$L$ is an affine automorphism such that $L(\{y=0\})=l$.

Put $\gamma (x)=H(x,0)$. By assumptions of the theorem the
mapping $\gamma :k\to k^2$ is an injection and
$\gamma '(x)\neq 0$ for $ x\in k$, hence $\gamma$ is an
embedding of the line in the plane. By Abhyankar-Moh theorem
\cite{1} there exist an automorphism $H_1\in \Aut$ such that
$\gamma (x)=H_1(x,0)$. Let $G=H_1\circ H$.
We get $\Jac G=\Jac H_1^{-1}\Jac H$ is a non-zero constant and
$G(x,0)=(x,0)$, so by lemma~2.3  $G\in \Aut$.
Therefore $H=H_1\circ G$ is a polynomial automorphism.
\enddemo

\head 		Injektywno\'s\'c na jednej prostej
\endhead        

\subhead Streszczenie
\endsubhead

{\eightpoint
Niech $k$ b\e{}dzie cia\l{}em algebraicznie domkni\e{}tym
charakterystyki zero.
Niech $H:k^2\to k^2$ b\e{}dzie odwzorowaniem wielomianowym
kt\'orego Jakobian $\Jac H$ jest sta\l{}\a{} r\'o\z{}n\a{} od zera.
W pracy dowodzimy, \z{}e je\'sli istnieje prosta $l \subset k^2$
na kt\'orej $H|_l:l\to k^2$ jest injekcj\a{},
to $H$ jest automorfizmem wielomianowym. }
\enddocument